\documentstyle[12pt,amscd]{amsart}

\newtheorem{theorem}{Theorem}[section]
\newtheorem{lemma}[theorem]{Lemma}
\newtheorem{corollary}[theorem]{Corollary}
\newtheorem{proposition}[theorem]{Proposition}

\def\ini{\operatorname{in}}
\def\char{\operatorname{char}}
\def\gin{\operatorname{gin}}
\def\Lex{\operatorname{Lex}}
\def\embdim{\operatorname{embdim}}
\def\depth{\operatorname{depth}}
\def\reg{\operatorname{reg}}
\def\Ext{\operatorname{Ext}}

\def\geom{\operatorname{g-reg}}

\def\To {\longrightarrow}
\def\mm{{\frak m}} 
\def\p{{\frak p}}

\begin{document}

\title{Castelnuovo-Mumford Regularity\\ and finiteness of Hilbert Functions}
\author{Maria Evelina Rossi, Ngo Viet Trung and Giuseppe Valla}
\address{Department of Mathematics,
University of Genoa, Via Dodecaneso 35, 16146 Genoa, Italy}
\email{rossim@@dima.unige.it}
\address{Institute of Mathematics, Vien Toan Hoc, 18 Hoang Quoc Viet, 10307 Hanoi, Vietnam}
\email{nvtrung@@math.ac.vn}
\address{Department of Mathematics,
University of Genoa, Via Dodecaneso 35, 16146 Genoa, Italy}
\email{valla@@dima.unige.it}
\thanks{The first and third authors are partially supported by MIUR of
Italy. The second author is partially supported
by the National Basic Research Program of Vietnam}

\maketitle

\section*{Introduction}

The Castelnuovo-Mumford regularity is a very important invariant of graded modules which arises naturally in the study of finite free resolutions. There have been several results which establish bounds for the Castelnuovo-Mumford regularity of projective schemes in terms of numerical characters. Unfortunately, these invariants  are often difficult to handle and the problem of finding good bounds in terms of simpler invariants is a topic featured in much recent research.  \par
 
The notion of regularity has been used by S. Kleiman in the construction of bounded families of ideals or sheaves with given Hilbert polynomial, a crucial point in the construction of Hilbert or Picard scheme. In a related direction, Kleiman proved that if $I$ is an equidimensional  reduced  ideal in a polynomial ring $S$ over an algebraically closed field, then the coefficients of the Hilbert polynomial of $R=S/I$ can be bounded by the dimension and the multiplicity of $R$ (see \cite[Corollary 6.11]{Kl}).  Srinivas and Trivedi proved that the corresponding result does not hold for a local domain. However, they proved that there exist a finite number of Hilbert functions for a local Cohen-Macaulay ring of given multiplicity and dimension (see \cite{ST2}).
 The proofs of the above results are very difficult and involve deep results from Algebraic Geometry.\par
 
The aim of this paper is to introduce a unified approach which gives more general results and  easier proofs of the above mentioned results. This approach is based on a classical result of Mumford on the behaviour of the geometric regularity by hyperplane sections in \cite{Mu} which allows us to  bound the regularity once there is a bound for the size of certain component of the first local cohomology modules or a uniform bound for the Hilbert polynomials. The finiteness of the Hilbert functions then follows from the boundness of the regularity. More precisely, we shall see that  a class $\mathcal C$ of standard graded algebras has a finite number of Hilbert functions if and only if there are upper bounds for the regularity and the embedding dimension of the members of $\mathcal C$. \par

In Section 1 we will prepare some preliminary facts on Castelnuovo-Mumford regularity and related notions such as weak regularity and geometric regularity. In Section 2 we will clarify the relationship between the boundedness of Castelnuovo-regularity (resp. geometric regularity) and the finiteness of the Hilbert functions (resp.  Hilbert polynomials). The main technique is the theory of Gr\"obner basis.\par

In Section 3 we will show that the boundedness of the regularity of classes of graded algebras with positive depth can be deduced from the boundedness of the dimension of  the zero-graded component of the first local cohomology modules. As a consequence, Kleiman's result follows from the well-known fact that for equidimensional reduced schemes, the dimension of the first sheaf cohomology module is equal to the number of the connected components minus one. We also give an example showing that Kleiman's result does not hold if the graded algebras are not equidimensional.\par

The local case is studied in Section 4. First we present a local version of Mumford's Theorem,  which allows us  to  reduce the problem of bounding the regularity of the tangent cone to the problem of bounding the Hilbert polynomial. Such a bound for the Hilbert polynomials was already established in \cite{RVV}, \cite{Tr}. From this we can easily deduce Srinivas and Trivedi's result on the finiteness of the Hilbert functions of Cohen-Macaulay local rings with given dimesion and multiplicity. This approach can be used to prove the more general result that the number of numerical functions which can arise as the Hilbert functions of local rings with given dimension and extended degree is finite (see \cite{RTV} for details).

\section{Variations on Castelnuovo-Mumford~regularity}

Throughout this paper let $S = k[x_1,\ldots,x_r]$ be a polynomial 
ring over a field $k.$ 

\noindent Let $M=\oplus_tM_t$ be a finitely generated graded $S$-module and let 
 $$0 \To  F_s \To  \cdots \To  F_1 \To  F_0 \To  M \To 0$$
be a minimal graded free resolution of $M$ as $S$-module.
Write  $b_i$ for the maximum of the  degrees of the generators of $F_i$.
Following \cite{E}, Section 20.5, we say that $M$ is {\it m-regular} for 
some integer $m$ if $b_j-j\le
m$ for all $j$.

The {\it Castelnuovo-Mumford regularity} $\reg(M)$ of $M$ is defined 
to be the least integer $m$ for which $M$ is $m$-regular, that is,
$$\reg(M) = \max\{b_i-i|\ i = 0,\ldots,s\}.$$

It is well known that $M$ is $m$-regular if and only if 
$\Ext_S^i(M,S)_n=0$ for all $i$ and all $n\le -m-i-1$ (see \cite{EG}).
This result is hard to apply because in principle infinitely many 
conditions must be checked. However, in some cases, it suffices to 
check just a few.  

We say that $M$ is {\it weakly m-regular} if $\Ext_S^i(M,S)_{-m-i-1}=0$ for
all $i\ge 0$.  Due to a result of Mumford, if $\depth M > 0$ and $M$
is weakly $m$-regular, then $M$ is $m$-regular (see \cite[20.18]{E}). From this we can easily deduce the following result (see \cite[Corollary 1.2]{RTV}).

\begin{proposition} \label{basic1} Let $R=S/I$ be a standard graded algebra and  $ m $ a non-negative
integer. If $R$ is weakly $m$-regular, then $R$ is $m$-regular.
\end{proposition}

Using local duality we can characterize these notions of regularity by means of the local cohomology modules of M.

Let $\Im$ denote the maximal graded ideal of $S$ and 
let $M$  be a finitely generated graded $S$-module.
For any integer $i$ we denote by $H_{\Im}^i(M)$ the $i$-th local 
cohomology module of $M$ with respect to $\Im$.
By local duality (see \cite{E},
A4.2) we have
$$H_{\Im}^i(M)_m \cong \Ext_S^{r-i}(M,S)_{-m-r}$$ for all $i$ and $m$.
Thus, $M$ is $m$-regular if and only if
$H_{\Im}^i(M)_n =0$ for all $i$ and $n \ge m-i+1$, and
$M$ is weakly $m$-regular if and only if $H_{\Im}^i(M)_{m-i+1}=0$ for all
$i$. In particular,
$\reg(M)$ is the least integer $m$ for which
$H_{\Im}^i(M)_n = 0$ for all $i$ and $n \ge m-i+1$. Hence the
Castelnuovo-Mumford regularity can be defined for any
finite $S$-module  regardless of its presentation. 

For any integer $i$ we set $a_i(M) := \max\{n|\ H_{\Im}^i(M)_n
\neq 0\},$ where $a_i(M) = -\infty$ if
$H_{\Im}^i(M) = 0$. Then
$$ \reg(M)  =  \max\{a_i(M)+i|\ i \ge 0\}.$$ 

A  relevant remark is that the Castelnuovo-Mumford regularity controls the behaviour of the Hilbert
function.
We recall that
the Hilbert function of $M$ is the numerical function $$h_M(t)=\dim_k(M_t).$$  
The Hilbert polynomial $p_M(X)$ of $M$ is the polynomial of degree $d-1$ such that  we have $h_M(t)=p_M(t) $  for $t>>0.$ 
 
The Hilbert function  $h_M(n)$
and the Hilbert polynomial  $p_M(n)$ are related by the 
formula $h_M(n) = p_M(n)$ for $n > \reg(M).$ This is a consequence of the Serre formula which holds for every integer $n$:
$$h_M(n) - p_M(n) = \sum_{i\ge 0} (-1)^i \dim_kH_{\Im}^i(M)_n.$$
For a proof, see for example \cite[Theorem 4.4.3]{BH}.\medskip

Inspired by the notion of regularity for sheaves on projective spaces, it is natural to introduce the following weaker notion  of regularity. 

 We say that $M$ is {\it geometrically
$m$-regular} if $H_{\Im}^{i}(M)_n=0$ for all $i>0$
and $n \ge m-i+1$, and we define the {\it geometric regularity} $\geom(M)$
of $M$ to be the least integer $m$ for
which $M$ is geometrically $m$-regular. 

It is clear that $$\geom(M) =
\max\{a_i(M)+i|\ i > 0\}.$$
Hence we  always have 
$$\geom(M) = \reg(M/H_{\Im}^0(M)) \le \reg(M).$$
In particular, $\reg(M) = \geom(M)$ if $\depth M >  0$.  

For a standard graded algebra $R=S/I, $ a theorem of Gotzmann  gives an upper bound  for  the geometric regularity in term of an integer  which can be computed  from   the Hilbert polynomial  of $R.$ 

\begin{theorem} \label{G} {\rm (see \cite{G})}
Assume that
 $$p_R(n)=\binom {n+a_1}{a_1}+\binom {n+a_2-1}{a_2}+\dots + \binom {n+a_s-(s-1)}{a_s}$$  with $a_1\ge a_2 \ge \dots \ge a_s \ge 0.$ Then   
$$\geom(R) = \reg(S/I^{sat})  \le s-1.$$
\end{theorem}

For example, if $R$ has dimension $1$ and multiplicity $e$, then its Hilbert polynomial is $$p_R(n)=e=\binom {n}{0}+\binom {n-1}{0}+\dots + \binom {n-(e-1)}{0}$$ so that $\geom(R)\le e-1.$ In particular, if $R$ is Cohen-Macaulay of dimension $1$ and multiplicity $e$, then $\reg(R)\le e-1.$\medskip

Unlike the regularity,  the geometric regularity does not behave well  under generic hyperplane sections.
Take for example $R=k[x,y,z]/(x^2,xy)$.
Then $\geom(R)=\reg(R)=1$ while $\geom(R/zR)=0,$ $ \reg(R/zR)=1.$

However the following result of Mumford (see \cite[page 101, Theorem]{Mu})
gives us the possibility to control this  behaviour. 
It  will be the basic result for our further investigation on the Castelnuovo-Mumford regularity of a standard graded algebra $R.$
A ring theoretic proof   can be found in \cite[Theorem 1.4]{RTV}.

\begin{theorem}\label{M}  Let $R=S/I$ be a standard graded algebra and $z\in R_1$ a regular linear form in $R.$ If $\geom(R/zR)\le m,$ then

{\rm (a)}  for every $s\ge m+1$,
$$\dim_k H_{\Im}^1(R)_m =\dim_k H_{\Im}^1(R)_s +\sum_{j=m+1}^s\dim_k H_{\Im}^0(R/zR)_j,$$

{\rm (b)} $\reg(R)\le m+\dim_k H_{\Im}^1(R)_m,$

{\rm (c)} $\dim_k H_{\Im}^1(R)_t =p_R(t)-h_R(t)$ for every $t\ge m-1$. 
\end{theorem}

\section{Finiteness of Hilbert Functions}  

The aim of this section is to clarify how the finiteness of Hilbert
functions for a given class of standard graded algebras is related to  the Castelnuovo-Mumford regularity of the members of the class.

In the following we let $R=S/I$ where $S=k[x_1,\ldots,x_r]$ is a polynomial ring over an infinite field $k$ of any characteristic and $I$ an homogeneous ideal of
$S.$ We will denote by $\embdim(R) $ its embedding dimension, that is $h_R(1).$ 

It is well-known that 
$$h_{S/I}(t) = h_{S/\ini(I)}(t)$$
for all $t \ge 0$ where $\ini(I)$ denotes the initial ideal of $I$ with respect to some term order. Therefore, we can pass to our study to  Hilbert function of factor rings by initial ideal.  
 
On the other hand, we have the following basic result of Bayer and Stillman on the behaviour of the regularity when passing to initial ideals.  

\begin{proposition}\label{S} {\rm (see \cite {BS}) } Let $R=S/I $ be a  standard graded algebra. Then
$$\reg(R) = \reg(S/\gin(I)), $$
where $\gin(I) $ is the generic initial ideal of $I$ with respect to 
 the reverse lexicographic order.
\end{proposition}

Moreover, it follows from a result of Bigatti and Hulett (for $\char(k) = 0$) and of Pardue  (for $\char(k) > 0$) that among all the ideals with the same Hilbert functions, the lex-segment ideal has the largest regularity.
Recall that the lex-segment ideal $\Lex(I)$ of $I$ is the monomial ideal which is generated in every degree $t$ by the first  $h_I(t)$ monomials in the lexicographical order.

\begin{proposition}\label{BH}  {\rm (see \cite{B}, \cite{Hul}, \cite{Pa})}
Let $R=S/I$ be a  standard graded algebra. Then
$$\reg(R) \le  \reg(S/\Lex(I)).$$
\end{proposition}
 
In the following  we will employ the following notations.

Let $\mathcal{C}$ be a class of standard graded algebras. We say:

\begin{itemize}
\item  $\mathcal{C}$ {\it is HF-finite} if the number of numerical functions which arise as the Hilbert functions of $ R \in \mathcal{C} $ is finite,
\item  $\mathcal{C}$ is {\it  HP-finite} if the number of  polynomials which arise as the Hilbert polynomials of $ R \in \mathcal{C} $, is finite,
\item $\mathcal{C}$ is {\it reg-bounded}  if  there exists an integer $t$ such that $\reg(R)\le t$ for all $R \in \mathcal{C}$, 
\item $\mathcal{C}$ is {\it  g-reg-bounded}  if  there exists an integer $t$ such that $\geom(R)\le t$ for all $R \in \mathcal{C}$, 
\item $\mathcal{C}$ is {\it  embdim-bounded}  if  if  there exists an integer $t$ such that $\embdim(R)\le t$ for all $R \in \mathcal{C}$.
\end{itemize}

Moreover,  we will denote by $\overline {\mathcal{C}}$ the class of graded algebras of the form $\bar R := R/H_{\Im}^0(R)$ for $R \in \mathcal{C}.$ Note that 
$$h_R(t) = h_{\bar R}(t)$$
for $t \gg 0$ and that 
$$\geom(R) = \geom(\bar R) = \reg(\bar R).$$
Then $\mathcal{C}$ is HP-finite if and only if $\overline {\mathcal{C}}$ is HP-finite and $\mathcal{C}$ is g-reg-bounded  if and only if $\overline{\mathcal{C}}$ is reg-bounded.

\begin{theorem}\label{BS} Let $\mathcal{C}$ be a class of standard graded algebras. Then \par 
{\rm  (a)} $\mathcal{C}$  is HF-finite if and only if $\mathcal{C}$   is reg-bounded and embdim-bounded,\par
{\rm (b)} $\mathcal{C}$ is HP-finite if and only if  $\mathcal{C}$  is g-reg-bounded and $\overline {\mathcal{C}}$ is embdim-bounded.
\end{theorem}

\begin{pf}  We may assume that the base field is infinite by tensoring it with a transcendental extension.\par

(a) If  $\mathcal{C}$ is HF-finite, then    $\mathcal{C}$   is 
embdim-bounded. By Proposition \ref{BH},  for every $S/I\in \mathcal{C}$ we have $\reg(S/I) \le  \reg(S/Lex(I)).$ Since there are a finite number of possible Hilbert functions for $S/I$,  there are also  a finite number of lex-segment ideals for $I$, so that $\mathcal{C}$ is reg-bounded by Proposition \ref{BH}.  

Conversely, assume $\mathcal{C}$   is reg-bounded and embdim-bounded. By Proposition \ref{S}, for all
$R=S/I \in \mathcal{C}$  we have
$$\reg(S/I)=\reg(S/\gin(I)) \ge D-1,$$ where $D$ is the maximum degree of the  monomials in the standard set of generators of $\gin(I).$ Since  the embedding dimension is bounded,   the number of monomials of degree smaller than or equal to
$D$ in $S$ are finite. Thus, there are only a finite number of possibilities for $\gin(I)$ and hence also for the Hilbert functions of $S/I$ because $h_{S/I}(t) = h_{S/\gin(I)}(t)$.
 
(b) By Theorem \ref{G} we have $\geom(S/I )\le s-1 $ where the integer $s$  depends only on the Hilbert polynomial of
$S/I.$   Therefore, if  $\mathcal{C}$ is HP-finite, then $  \mathcal{C}$ is g-reg-bounded. As remarked above, $\overline{\mathcal{C}}$ is  reg-bounded. 
Furthermore, for every $R \in \overline{ \mathcal{C}}$  we have
$h_R(1) \le h_R(n) = p_R(n)$ for $n = \reg(R).$  Since  $\overline{\mathcal{C}}$ is HP-finite, there are only a finite number of Hilbert polynomials $p_R(n)$. Let $t = \max\{p_R(n)| n = \reg(R), R \in \overline{ \mathcal{C}}\}$. Then $h_R(1) \le t$ for all $R \in \overline{ \mathcal{C}}$. Hence $\overline{\mathcal{C}}$ is embdim-bounded. 

Conversely, if $\mathcal{C}$ is g-reg-bounded, then $\overline{ \mathcal{C}}$ is reg-bounded. If moreover
$\overline {\mathcal{C}}$ is embdim-bounded, then $\overline {\mathcal{C}}$ is HF-finite by (a). This implies that $ \mathcal{C}$ is HP-finite.
\end{pf}

\begin{corollary}  Let $\mathcal{C}$ be a class of standard graded algebras.
 Then $\mathcal{C}$  is HP-finite if and only if    $\overline{ \mathcal{C}}$   is HF-finite.
\end{corollary}
 
The following example shows that  $\overline{\mathcal{C} } $ is embdim-bounded does not imply  $\mathcal{C}  $ is
embdim-bounded.

Let  $\mathcal{C}  $  be the class of algebras  of the form
$$ R_n = k[x_1,...,x_n]/(x_1^2,...,x_{n-1}^2,x_1x_n,...,x_{n-1}x_n) $$
for $n > 0$. Then $ \overline R_n  \cong k[x_n].$
Hence $\overline {\mathcal{C}}  $ is embdim-bounded while  $\mathcal{C}$ is not because    $\embdim R_n  = n.$  
 
\section{Reg-bounded algebras and Kleiman's~ Theorem.}

 The aim of this section is to present a relevant class of algebras  which are
reg-bounded and HF-finite. As an application, 
  we give a  proof of a   theorem of Kleiman (see \cite[Corollary 6.11]{Kl}),  which says 
that the class  of   graded reduced  and equidimensional algebras  with given
multiplicity and dimension is reg-bounded and HF-finite.  
The main tool is the afore mentioned result of Mumford (Theorem \ref{M}).

For every integer $p \ge 1$ we define recursively the following polynomials $F_p(X)$ with rational coefficients. We let 
$$F_1(X) := X $$ 
and, if $p\ge 2,$ then we let 
$$F_p(X) :=F_{p-1}(X)+X\binom{F_{p-1}(X)+p-1}{p-1}.$$
 
\begin{theorem} \label{T1}  Let $\mathcal{C}$ be a class of standard graded
algebras with the following properties:

\begin{enumerate}
\item $\depth(R) > 0$ for every  $R \in  \mathcal{C}$,
\item for every  $R \in  \mathcal{C}$ with $\dim R\ge 2, $ there exists  some  regular linear form $x \in R$ such that  $R/(x)^{sat} \in  \mathcal{C}$, 
\item there exists an integer $t$ such that   $\dim_k H_{\Im}^1(R)_0 \le t$ 
for every  $R \in  \mathcal{C}.$
\end{enumerate}

\noindent Then for every  $R \in \mathcal{C} $ with $d = \dim R$ we have
$$\reg(R)\le F_{d }(t).$$ 
\end{theorem}

 It would be interesting the problem of finding the complexity of the regularity bound. Explicit bounds  can be found in \cite{RTV}.

We need the following observation for the proof of the above theorem.
 
\begin{lemma}\label{P} Let $\mathcal{C}$ be a class of standard graded algebras 
as in the above theorem.
Then for all $j \ge 0 $ and for all $R \in \mathcal{C}$ we have
$$\dim_k H_{\Im}^1(R)_j \le t \binom{j+d-1}{d -1} .$$
\end{lemma}

\begin{pf} If $j = 0$, the conclusion holds by the assumption. If $d = 1$, then $R$ is a one-dimensional  standard graded algebra of positive depth. Hence $H_{\Im}^0(R)=0$ so
that by Serre formula we get for every $ j $,
$$\dim_k H_{\Im}^1(R)_j =p_R(j)-h_R(j)=e-h_R(j), $$ 
where $e$ denotes the multiplicity of $R.$ In particular, if $R \in \mathcal{C} $ and $j\ge 0, $ we have 
$$\dim_k H_{\Im}^1(R)_j =e-h_R( j) \le e-1 =\dim_k H_{\Im}^1(R)_0  \le t.$$
Hence  the conclusion holds for every $R \in \mathcal{C} $ of dimension $1$.  \par

Now let $j\ge 1$ and $R \in \mathcal{C} $ of dimension $d\ge 2$.  We may identify $R$ with the graded flat extension $R'$ of Theorem \ref{T1} (2). Then there exists  a regular element $x\in R_1$ such that $ B=R/(x)^{sat} \in \mathcal{C}.$  It is clear that $H^i(R/xR)=H^i(B)$ for every
$i>0.$ From the short exact sequence 
$$0 \To  R(-1) \overset x \To  R \To  R/xR \To  0$$
we get for every $j$ an exact sequence $$\cdots \to H_{\Im}^1(R)_{j-1}\to H_{\Im}^1(R)_{j}\to H_{\Im}^1(R/xR)_j\to \cdots$$
Hence $$\dim_k H_{\Im}^1(R)_j \le
 \dim_k H_{\Im}^1(R)_{j-1} +\dim_k H_{\Im}^1(R/xR)_j =\dim_k H_{\Im}^1(R)_{j-1} +\dim_k H_{\Im}^1(B)_j.$$
Since $\dim B = d-1, $  we can use the inductive assumption and we get
$$\dim_k H_{\Im}^1(R)_j \le t\binom{j+d-2}{d-1}+t\binom{j+d-2}{d-2}=t\binom{j+d-1}{d-1}.$$
\end{pf}

\begin{pf}{\it{of Theorem \ref{T1}.}}  If $R \in \mathcal{C}$  is  one-dimensional, then $R$ is Cohen-Macaulay. Using Theorem \ref{G} and  the proof of Lemma \ref{P} we have
$$\reg(R) = \geom(R) \le e-1=\dim_k H_{\Im}^1(R)_0 \le t=F_1(t).$$  
If $R \in \mathcal{C} $ is of dimension $d\ge 2$, we replace $R$ by the graded flat extension $R'$ of (2).  Then there exists
a regular element $x\in R_1$ such that
$B=R/(x)^{sat} \in \mathcal{C}.$ We have $\depth(B)\ge 1$ and if we let 
$m:=\reg(B),$ then 
$$m  = \geom(B) = \geom(R/xR).$$ 
By Theorem \ref{M} and Lemma \ref{P}, this implies 
$$\reg(R)\le m+\dim_k H_{\Im}^1(R)_m 
\le m+t\binom {m+d-1}{d-1}.$$
Since $\dim B=d-1,$  by induction we have $m \le F_{d-1}(t).$ Thus,  
$$\reg(R)\le F_{d-1}(t)+t\binom{F_{d-1}(t)+d-1}{d-1}=F_d(t).$$
\end{pf}

We want to apply now the above result to the class $\mathcal{C}$ of reduced equidimensional graded algebras with given multiplicity $e$ and dimension $d\ge 1$. 

It is clear that every $R \in  \mathcal{C}$ has positive depth. Moreover, we have the following easy lemma.

 \begin{lemma} \label{A} If $R$ is a reduced
equidimensional graded algebra  over an algebraically closed field with multiplicity $e$, then   
$\dim_k H_{\Im}^1(R)_0 \le e-1.$  \end{lemma}
 
 \begin{pf} If $\dim R = 1$, then $R$ is Cohen-Macaulay and, as above, $\dim_k H_{\Im}^1(R)_0 =e-1. $ 
If $\dim R \ge 2$, then  
$\dim_k H_{\Im}^1(R)_0 = N-1$ where $N$ is the number of the connected components of the corresponding scheme \cite[Theorem 1.2.6 (b)]{M}. Since $N \le e$, we get $\dim_k H_{\Im}^1(R)_0 \le e -1.$
 \end{pf}

We can prove  now that  Kleiman's result is a  particular case of our theorem.

\begin{theorem} \cite[Corollary 6.11]{Kl}
Let $\mathcal{C}$ be the class of reduced equidimensional graded algebras over an algebraically closed field with given multiplicity $e$ and dimension $ \le d$, $d \ge 1.$  Then $\mathcal{C}$ is HF-finite. 
\end{theorem}

\begin{pf} By Bertini theorem \cite[Corollary 3.4.14]{FOV}, if $\dim R \ge 2$, there exists a regular linear form $x \in R$ such that $R/(x)^{sat}$ is  a reduced equidimensional algebra with $\dim R/(x)^{sat} = \dim  R -1$. Therefore, we may apply Theorem \ref{T1} and Lemma \ref{A} to get $ \reg(R) \le F_d(e-1)$ for every $R$ in the class $\mathcal{C}.$ So $\mathcal{C}$ is
reg-bounded. Now we need only to prove that  $\mathcal{C}$ is embdim-bounded since by Theorem \ref{S} (a), these conditions will imply  that $\mathcal{C}$ is finite. 

Let $R \in \mathcal{C}$ be arbitrary. Write $R = S/I$ where $S$ is a polynomial ring over an algebraically closed field and $I=\bigcap _{i=1}^r \p_i  $ is an intersection of equidimensional ideals $\p_i$. It is obvious that $\embdim R  \le \sum_{i=1}^r \embdim S/\p_i$ and $r \le e$. On the other hand, we know that $\embdim S/\p_i \le e_i+d-1$, 
where $e_i$ is the multiplicity of $S/\p_i$. Therefore,
$$\embdim R  \le \sum_{i=1}^r  e_i +d-1 =  e+r(d-1) \le ed.$$ 
\end{pf}

Theorem \ref{T1} does not hold if we delete the assumption that every element of the class is reduced. Take for example the class $\mathcal{C}$ of the graded algebras 
$$R_n :=k[x,y,z,t]/(y^2,xy,x^2,xz^n-yt^n)\; (n \ge 1).$$  
Note that $(y^2,xy,x^2,xz^n-yt^n)$ is a primary ideal.
We have $\dim (R_n)=2$ and $e(R_r)=2$. The minimal free resolution of $R_r$ over $S = k[x,y,z,t]$ is given by
$$0\to S(-n-3) \to S(-3)^2\oplus S(-n-2)^2\to S(-2)^3\oplus S(-n-1)\to S \to R_n \to 0.$$   Hence $\reg(R_n)=n$. Therefore, $\mathcal{C}$ is not reg-bounded and hence not HF-finite. \medskip

Theorem \ref{T1} does not hold  if we consider reduced graded algebras which are not necessarily equidimensional.
Let us consider the class of standard graded algebras $$R_n :=k[x,y,z_1,\dots,z_r]/(x) \cap (y,f_n) $$ where $f_n \in k[z_1,\dots,z_r] $ is an irreducible form of degree $n.$ 
We have  $\dim R_n = r+1,$  $e(R_n)=1 $, but $ \reg(R_n)=n.$

 \section{The local version of Mumford Theorem}  

Let $(A,\mm)$ be a local ring of dimension $d$ and multiplicity $e$. Let $$G=gr_{\mm}(A)=\oplus_{n\ge 0}(\mm ^n/\mm^{n+1})$$ be the associated graded
 ring of $A$. The Hilbert function and the Hilbert polynomial of $A$ is by definition the Hilbert function and the Hilbert polynomial of the standard graded algebra $G$, namely $h_A(t) :=h_G(t)=\lambda(\mm^t/\mm^{t+1})$ and $p_A(t) :=p_G(t)$.

We will need also the first iterated of these functions and polynomials. So we let $$h^1_A(t) :=h^1_G(t)=\sum_{j=0}^th_G(j)=\lambda(A/\mm^{t+1})$$ and we denote by
$p_G^1(t) $ the corresponding polynomial, that is the polynomial which verifies the equality  $h_G^1(t)=p_G^1(t) $ for
$t \gg 0. $  

Let $x$ be a superficial element in $\mm$ and let 
$$\bar G :=gr_{\mm/xA}(A/xA)$$ 
be the associated graded ring of the local ring $A/xA.$ Let $x^*=\overline{x}\in (\mm/\mm^2)$ be the initial form of $x$ in $G.$ We can consider the standard graded algebra $G/x^*G$ and compare it with 
$\overline{G}$. These two algebras are not the same, unless $x^*$ is a regular element in $G$, but they have the 
same geometric regularity, namely 
$$\geom(G/x^*G)=\geom(\bar G).$$ 
This has been proved in \cite[Lemma 2.2]{RTV}.

If $A$ has positive depth, then it is well known that $x$ is a regular element in $A$ and furthermore $$p_G(t)=p^1_{\overline{G}}(t).$$
In that case, one can prove (see \cite[Lemma 3.2]{RTV})  that $\reg(G)=\geom(G)$
(even when $G$ does not  necessarily have positive depth).

With the above notations we present now a local version of Mumford's Theorem \ref{M}.

\begin{theorem}\label{HF} Let $(A,\mm)$ be a local  ring  of dimension $d\ge 2$ and positive depth. Let $x$ be a superficial  element in $\mm$ and $m:=\reg(\overline{G})$. Then
$$\reg(G)\le m+h^1_{\overline{G}}(m).$$
\end{theorem}

\begin{pf} We have 
$$m=\reg(\overline{G})\ge \geom(\overline{G})=\geom(G/x^*G)=\geom(G/(H_{G_+}^0(G)+x^*G)),$$ 
where the last equality follows because
$(H_{G_+}^0(G)+x^*G)/x^*G)$ has finite length in $G/x^*G.$
Now we remark that $x^*$ is a regular element in $G/H_{G_+}^0(G),$ hence, 
by using Mumford theorem, we get
\begin{align*}
\reg(G) =\geom(G) & =\geom(G/(H_{G_+}^0(G)) \le m+p_{G/H_{G_+}^0(G)}(m)\\
& = m+p_G(m) =m+p^1_{\overline{G}}(m)=m+h^1_{\overline{G}}(m)
\end{align*}
where the equality
$p^1_{\overline{G}}(m)=h^1_{\overline{G}}(m)$ is a consequence of the fact that 
$m=\reg(\overline{G}).$
\end{pf}

Using this local version of Mumford's theorem we can easily deduce the result of  Srinivas and Trivedi which says that the number of Hilbert functions of Cohen-Macaulay local rings with fixed dimension and multiplicity is finite (see \cite{ST2}). For that we need  the following inequality proved in
\cite{RVV} and \cite{Tr}.

\begin{proposition} \label{in} 
Let $(A,\mm)$ be a local ring of dimension $d\ge 1$ and $J$ an ideal generated by a system of parameters in $\mm.$ Then 
$$H_A(n)\le \ell(A/J) \binom{n+d-2}{d-1}+\binom{n+d-2}{d-2}.$$ 
\end{proposition}

If  $A$ is  a $d$-dimensional Cohen-Macaulay ring of multiplicity $e$, then from the above Proposition we immediately get the
inequality 
$$ H_A(n)\le e\binom{n+d-2}{d-1}+\binom{n+d-2}{d-2}.$$

For every $d\ge 1$ we define recursively the following polynomials $Q_d(X)$ with rational coefficients. We let 
$$Q_1(X):=X-1$$ and, if $d\ge 2,$ then we let $$Q_d(X):=Q_{d-1}(X)+X\binom{Q_{d-1}(X)+d-2}{d-1}+\binom{Q_{d-1}(X)+d-2}{d-2}.$$

\begin{theorem} Let $(A,\mm)$ be a Cohen-Macaulay local ring of dimension $d\ge 1$ and multiplicity $e.$ Then $$\reg(G)\le Q_d(e).$$
\end{theorem}

\begin{pf} If $d=1$, then $\reg(G)\le e-1=Q_1(e).$ Let $d\ge 2$ and $x$ be a superficial element in $\mm$. Then $A/xA$ is a Cohen-Macaulay local ring of dimension $d-1$ and multiplicity $e$. By
Proposition \ref{in}  we get  $$H^1_{\overline{G}}(n)\le
 e \binom{n+d-2}{d-1}+\binom{n+d-2}{d-2}.$$ By Theorem \ref{HF} it follows that,  if $m=\reg(\overline{G}),$  then $$\reg(G)\le m+h^1_{\overline{G}}(m)\le
m+e\binom{m+d-2}{d-1}+\binom{m+d-2}{d-2}.$$ By induction we have $m\le Q_{d-1}(e),$ so that $$\reg(G)\le
Q_{d-1}(e)+e\binom{Q_{d-1}(e)+d-2}{d-1}+\binom{Q_{d-1}(e)+d-2}{d-2}=Q_d(e).$$
\end{pf}

\begin {corollary} The number of numerical functions which can arise as the Hilbert functions of Cohen-Macaulay local rings with given dimension and multiplicity is finite.
\end{corollary}

\begin{pf} By a classical result of Abhyankar we have $v(\mm) \le e+d-1$, where $v(\mm)$ is the embedding dimension of $G$.  Now we need only to apply  Proposition \ref{BS} and the above theorem. \end{pf}

The analogous of Kleiman result does not hold in the local case. Srinivas and Trivedi gave the following example showing that classes of local domains of fixed dimension and multiplicity need not to be HF-finite.

Let $$A_r:=k[[x,y,z,t]]/(z^rt^r-xy, x^3-z^{2r}y,y^3-t^{2r}x,x^2t^r-y^2z^r).$$ 
It is easy to see that
$A_r$ is a local domain and the associated graded ring of $A_r$ is the standard graded algebra 
$$G_r=k[x,y,z,t]/(xy,x^3,y^3,x^2t^r-y^2z^r).$$
We have $\reg(G_r)=r+1$ and
$$H_{A_r}(n) = \left\{\begin{array}{lll} 5n-1 & \text{for} & n \le r,\\
4 n + r  & \text{for} & n > r. \end{array} \right.$$

Finally, we remark that the above approach can be used to prove  that the number of numerical functions which can arise as the Hilbert functions of local rings with given dimension and extended degree is finite. Note that  extended degree coincide with the usual multiplicity for Cohen-Macaulay local rings. We refer to \cite{RTV} for details. Furthermore, one can prove similar results for  Hilbert functions of finitely generated modules over local rings with respect to $\mm$-primary ideals (see \cite{T1, T2, Li}).

\end{document}